**Oscillations in Mertens Theorems and Other Finite Sums and Products**
**N. A. Carella, July, 2013.**

*Abstract:* This note simplifies the proof of a recent result on the oscillation of the prime product in Mertens Theorem, and provides a quantitative expression for the error term. In addition, the corresponding oscillations results for the finite sum of the reciprocals of the prime numbers up to a fixed number and other finite sums and products are also given.

Mathematics Subject Classifications 2000: 11A25, 11L20, 11N05.

Keywords: Prime Numbers, Arithmetic Functions, Sums over the Primes, Oscillation Theorems.

**1. Introduction**

For small real numbers $x \geq 1$, the numerical evaluations of the prime product $P(x) = (\log x)^{-1} \prod_{p \leq x} (1 - 1/p)^{-1}$ approach the real number $e^\gamma$ from above. The function $P(x)$ oscillates infinitely often symmetrically about the real number $e^\gamma = 1.781072417987...$ as $x$ tends to infinity. The number $e^\gamma$ is the limit of the product $P(x)$ as $x \to \infty$, see Theorem 7. Recently the following oscillation theorem was proved.

*Theorem 1.* ([DP]) The quantity $x^{1/2} \left( \prod_{p \leq x} (1 - 1/p)^{-1} - e^\gamma \log x \right)$ attains arbitrary large positive and negative values as $x \to \infty$.

This note significantly simplifies the proof of Theorem 1, and provides a quantitative expression for the error term; Theorem 1 does not address the quantitative aspect of the oscillations, see Section 4. In addition, the oscillations results for the corresponding finite sum of the reciprocals of the prime numbers $p \leq x$ and other finite sums and products are also given, see Sections 3 and 5. More precisely, the following version of Mertens Theorem, and other similar results will be confirmed.

*Theorem 2.* Let $x \geq x_0$. Then

$$\sum_{p \leq x} \frac{1}{p} = \begin{cases} \log\log x + B_1 + O(e^{-.2098(\log x)^{3/5} (\log\log x)^{-1/5}}) & \text{if } x \notin E, \\ \log\log x + B_1 + \Omega_\pm(x^{-1/2} \log\log\log x / \log x) & \text{if } x \in E, \end{cases} \quad (1)$$

where $B_1 > 0$ is a constant, and $E \subset \mathbb{R}$ is an exceptional subset of real numbers, see Definition 5. In particular, the expression

$$x^{1/2} \log x \left( \sum_{p \leq x} 1/p - \log\log x - B_1 \right) = \Omega_\pm(\log\log\log x) \quad (2)$$

attains arbitrary large positive and negative values as $x \to \infty$.

Its complementary result appears in Theorem 11. The earliest result on the oscillatory nature of the finite sums

$$\sum_{p^n \leq x} \frac{1}{p^n} \quad \text{and} \quad \sum_{p^n \leq x} \frac{\log p}{p^n}, \tag{3}$$

where $n \geq 1$ is a fixed parameter, appears to be the work of Phragmen, and probably other authors too, see [NW, p. 122]. The recent quantitative improvements for the finite sums are given in [RN], and for the product in [DP].

Intuitively and theoretically speaking, the oscillations are caused by both the excess of short prime gaps and the excess of large prime gaps in the intervals $[1, x]$ as $x \to \infty$. To confirm this claim, observe that replacing the sequence of prime numbers $p_1, p_2, \ldots, p_n \leq x$ with another sequence of numbers $a_1, a_2, \ldots, a_n \leq x$ of the same density as the sequence of prime numbers in the set of nonnegative numbers $\mathbb{N}$, but shorter gaps $a_{n+1} - a_n \leq p_{n+1} - p_n$, and then evaluating both products yield

$$P(x) \leq P_a(x) = (\log x)^{-1} \prod_{a \leq x} (1 - 1/a)^{-1}. \tag{4}$$

This suggests that if the prime numbers $p_1, p_2, \ldots, p_n \leq x$ in the interval $[1, x]$ have sufficiently many short prime gaps $p_{n+1} - p_n$, then $\pi(x) - li(x) \geq 0$. In contrast, replacing the sequence of prime numbers $p_1, p_2, \ldots, p_n \leq x$ with another sequence of numbers $b_1, b_2, \ldots, b_n \leq x$ of the same density, but larger gaps $p_{n+1} - p_n \leq b_{n+1} - b_n$, and then evaluating the products yield

$$P_b(x) = (\log x)^{-1} \prod_{b \leq x} (1 - 1/b)^{-1} \leq P(x). \tag{5}$$

Similarly, this suggests that if the prime numbers $p_1, p_2, \ldots, p_n \leq x$ in the interval $[1, x]$ have sufficiently many large prime gaps $p_{n+1} - p_n$, then $\pi(x) - li(x) < 0$. Advanced materials on the randomness of the prime numbers appear in [GV] and similar references.

The oscillation phenomenon can be simulated in numerical experiments. Id est, a composite sequence of numbers $c_1, c_2, \ldots, c_n \leq x$ of the same density as the sequence of prime numbers, but with pseudorandom strings of small gaps $c_{n+1} - c_n < p_{n+1} - p_n$ and large gaps $p_{n+1} - p_n < c_{n+1} - c_n$, both in the range 2 to $O(\log^2 c_n)$, and $|c_n - p_n| = O(\log^2 c_n)$ as specified by Cramer conjecture, can simulates the oscillations of the sequence of prime numbers in intervals of practical sizes, for instance, $x = 10^{10}$.

## 2. Auxiliary Materials

Let $\pi(x) = \#\{ p \leq x : p \text{ is prime} \}$ denotes the prime counting function. The current version of the Prime Number Theorem claims the following.

***Theorem* 3.** (Korobov, Vinogradov) Let $x \geq x_0$ be a large real number. Then

$$\pi(x) = li(x) + O(xe^{-.2098(\log x)^{3/5}(\log\log x)^{-1/5}}), \tag{6}$$

where $li(x)$ denotes the logarithm integral.

Prior to this result, there was the delaVallee Poussin form $\pi(x) = li(x) + O(xe^{-c(\log x)^{1/2}})$, $c > 0$ constant, and other earlier versions such as the Legendre form $\pi(x) \approx x/(\log x - B)$, $B > 0$ constant, and the Gauss form $\pi(x) \approx x/\log x$. The Riemann form $\pi(x) = li(x) + O(x^{1/2} \log x)$ provides an optimal asymptotic formula, confer [NW], [IV] and [FD] for details and proofs.

The oscillations of the function $\pi(x)$ have been investigated for many years, confer [RL], [HN], [KS], [SD] and other sources for details. The prototypical result on the difference $\pi(x) - li(x)$ accounts for the absolute minimum variations for infinitely many real numbers $x > 0$.

***Theorem* 4.** (Littlewood) (i) For large $x > 0$, the error term $\pi(x) - x/\log x = \Omega(x/(\log x)^2)$ is strictly positive.
(ii) The error term $\pi(x) - li(x) = \Omega_\pm(x^{1/2} \log\log\log x / \log x)$ is oscillatory.

A full proof is given in [EL, p. 191] and a weaker result appears in [IV, p. 306], see also [MV], et cetera. The actual intervals of sign changes are given in [SD] and earlier references, and a lower estimate $V(x) \geq c\log x$, $c > 0$ constant, for the number of sign changes per interval $[1, x]$ is given in [KS]. The numerical data for the difference $\pi(x) - li(x)$ up to the current limit of computing technology are given in [KT]. The lower estimate for the number $V(x)$ of sign changes implies that the equation $\pi(x) - li(x) = 0$ has infinitely many zeros. Additional advanced materials on the general theory of oscillations of function appear in [GR], [KS], et al.

***Definition* 5.** The exceptional subset of real numbers $E \subset \mathbb{R}$ is defined by $E = \{ x \in \mathbb{R} : |\pi(x) - li(x)| \geq x^{1/2} \log\log\log x / \log x \}$.

Other than the qualitative description of the size of the exceptional subset $E \subset \mathbb{R}$ as being an infinite subset, there appears to be no quantitative information in the literature on the precise nature of this subset of numbers. The subsets of numbers

$$E_+(x) = \{ x \in \mathbb{R} : \pi(x) - li(x) \geq x^{1/2} \log\log\log x / \log x \} \tag{7}$$

and

$$E_-(x) = \{ x \in \mathbb{R} : \pi(x) - li(x) \leq -x^{1/2} \log\log\log x / \log x \}$$

have different densities, see [BS]. Presumably, the subset $E$ is uncountable, but the intersection $E \cap \mathbb{N}$ might have zero density in the set of integers $\mathbb{N} = \{ 0, 1, 2, 3, \ldots \}$. The information encapsulated in Corollary 8 seems to imply that the subset of integers $E \cap \mathbb{N}$ is a very thin subset of $\mathbb{N}$.

The corresponding results for the differences $\pi(x, a, q) - li(x)/\varphi(q)$ for primes in arithmetic progressions are quite similar to the special case $a = 1$, $q = 1$, but more complicated and depend on the residue classes $a \bmod q$, and other distribution properties of the prime numbers, see [KS, Theorem 1.2], [FR], and [BS] for some details.

By virtues of Theorems 3 and 4, the unconditional approximation of the prime counting measure

Oscillations in Finite Sums and Products$$d\pi(x) = \begin{cases} d(li(x) + O(xe^{-.2098(\log x)^{3/5}(\log\log x)^{-1/5}})) & \text{if } x \notin E, \\ d(li(x) + \Omega_{\pm}(x^{1/2}\log\log\log x/\log x)) & \text{if } x \in E, \end{cases} \qquad (8)$$

holds for all large numbers $x \in \mathbb{R}$. And assuming the Riemann Hypothesis, the best possible approximation is

$$d\pi(x) = \begin{cases} d(li(x) + O(x^{1/2}\log x)) & \text{if } x \notin E, \\ d(li(x) + \Omega_{\pm}(x^{1/2}\log\log\log x/\log x)) & \text{if } x \in E. \end{cases} \qquad (9)$$

### 3. Sums and Products over the Primes

Let $\mathbb{P} = \{\, p \in \mathbb{N} : p \text{ is prime }\} \subset \mathbb{N}$ be the subset of prime numbers, and let $\mathbb{C} = \{\, s = \sigma + it : \sigma, t \in \mathbb{R}\,\}$ be the set of complex numbers. The tasks of determining the oscillations of a function $f : \mathbb{P} \to \mathbb{C}$ is often a long, winding, and difficult task. However, assuming various results and the Prime Number Theorem, the scheme involving the prime counting measure in conjunction with the Stieltjes integral provide a simple and unified approach to the determinations of the oscillation estimates of a large class of sums and products over the primes. Furthermore, any advances in the Prime Number Theorem immediately cascade into sharper estimates for these sums and products over the primes.

The verification of Theorem 2 is derived from Theorems 3 and 4. Subsequently, this is used to obtain the oscillation results for the corresponding products.

***Proof of Theorem* 2:** In terms of the prime counting measure $d\pi(t)$, the finite sum has the integral representation

$$\sum_{p \leq x} \frac{1}{p} = \sum_{n \leq x} \frac{\pi(n) - \pi(n-1)}{n} = \int_2^x \frac{d\pi(t)}{t}, \qquad (10)$$

see [RS, p. 67] for similar calculations using Stieltjes integral. To simplify the notations, only the exceptional part of the measure $d\pi(t)$ will be given, the calculations for the nonexceptional part is similar. Using (8) and integration by parts yield

$$\begin{aligned}\sum_{p \leq x} \frac{1}{p} &= \int_2^x \frac{d(li(t) + \Omega_{\pm}(t^{1/2}\log\log\log t/\log t))}{t} \\ &= \int_2^x \frac{dt}{t\log t} + \int_2^x \frac{d\Omega_{\pm}(t^{1/2}\log\log\log t/\log t)}{t} \\ &= \log\log x + c_1 + \frac{\Omega_{\pm}(x^{1/2}\log\log\log x/\log x)}{x} - \int_2^x \frac{\Omega_{\pm}(t^{1/2}\log\log\log t/\log t)dt}{t^2}, \end{aligned} \qquad (11)$$

where $c_1 = -\log\log 2$, $c_2 > 0$ are constants. Recall that the logarithmic integral $li(x) = \int_2^x dt/\log t$, $x \in \mathbb{R}$. The last integral is written in the more convenient form

$$\int_2^x \frac{\Omega_\pm(t^{1/2}\log\log\log t/\log t)dt}{t^2} = \int_2^\infty \frac{\Omega_\pm(t^{1/2}\log\log\log t/\log t)dt}{t^2} - \int_x^\infty \frac{\Omega_\pm(t^{1/2}\log\log\log t/\log t)dt}{t^2} \quad (12)$$
$$= c_2 + \Omega_\pm(x^{-1/2}\log\log\log x/\log x).$$

Now regroup all the terms and put in $B_1 = c_1 + c_2$ to establish the claim. ∎

A handful of different proofs of the standard version of Mertens Theorem for the finite sum

$$\sum_{p \leq x} 1/p = \log\log x + B_1 + O(1/\log x) \quad (13)$$

are given in the literature, see [HW], [VL], [DU], [TB] and related sources.

Basically, the asymptotic formula $\sum_{p \leq x} 1/p = \log\log x + B_1 + O(e^{-c(\log x)^{3/5}(\log\log x)^{-1/5}})$ was obtained by Vinogradov, see [VV]. Moreover, assuming the Riemann Hypothesis, the best possible asymptotic formula is $\sum_{p \leq x} 1/p = \log\log x + B_1 + O(x^{-1/2}\log x)$, see [RS]. An advanced algorithm for computing the prime harmonic sum $\sum_{p \leq x} 1/p$ is developed in [BH].

There are various analytical formulas for the Meissel–Mertens constant

$$B_1 = \gamma + \sum_p \log(1-1/p) + 1/p = 0.2614972128..., \quad (14)$$

and the Euler–Mascheroni constant

$$\gamma = \lim_{n \to x} \sum_{n \leq x} 1/n - \log x = 0.5772156649..., \quad (15)$$

see [HW], [FH], [LA], and [VL]. The analytic formulas are useful in analysis.

The extension of Mertens Theorem to primes in arithmetic progressions $qn + a$ with $\gcd(a, q) = 1$, $n \in \mathbb{N}$, has the form

$$\sum_{p \leq x,\, p \equiv a \bmod q} \frac{1}{p} = \frac{1}{\varphi(q)}\log\log x + B_{a,q} + O_q(1/\log x), \quad (16)$$

where $B_{a,q}$ is a constant, see [LC] and related papers.

***Theorem* 6.** Let $x \in \mathbb{R}$ be a large real number, then

(i) $\displaystyle\prod_{p \leq x}\left(1 - \frac{1}{p}\right)^{-1} = \begin{cases} e^\gamma \log x + O(e^{-.2098(\log x)^{3/5}(\log\log x)^{-1/5}}) & \text{if } x \notin E, \\ e^\gamma \log x + \Omega_\pm(x^{-1/2}\log\log\log x/\log x) & \text{if } x \in E, \end{cases}$ \quad (17)

(ii) $\displaystyle\prod_{p \leq x}\left(1 - \frac{1}{p}\right) = \begin{cases} (e^\gamma \log x)^{-1} + O(e^{-.2098(\log x)^{3/5}(\log\log x)^{-1/5}}) & \text{if } x \notin E, \\ (e^\gamma \log x)^{-1} + \Omega_\pm(x^{-1/2} \log\log\log x / \log x) & \text{if } x \in E, \end{cases}$

(iii) $\displaystyle\prod_{p \leq x}\left(1 + \frac{1}{p}\right)^{-1} = \begin{cases} (6e^\gamma \pi^{-2} \log x)^{-1} + O(e^{-.2098(\log x)^{3/5}(\log\log x)^{-1/5}}) & \text{if } x \notin E, \\ (6e^\gamma \pi^{-2} \log x)^{-1} + \Omega_\pm(x^{-1/2} \log\log\log x / \log x) & \text{if } x \in E, \end{cases}$

(iv) $\displaystyle\prod_{p \leq x}\left(1 + \frac{1}{p}\right) = \begin{cases} 6e^\gamma \pi^{-2} \log x + O(e^{-.2098(\log x)^{3/5}(\log\log x)^{-1/5}}) & \text{if } x \notin E, \\ 6e^\gamma \pi^{-2} \log x + \Omega_\pm(x^{-1/2} \log\log\log x / \log x) & \text{if } x \in E, \end{cases}$

as $x \to \infty$.

**Proof of (i):** Use the power series $-\log(1-z) = \sum_{n \geq 1} z^n/n$, $|z| < 1$, to expand the logarithm of the product over the prime $p \leq x$ into a power series. This step yields the expression

$$\sum_{p \leq x} \log(1 - 1/p)^{-1} = \sum_{p \leq x} \sum_{n=1}^\infty \frac{1}{np^n}$$
$$= \sum_{p \leq x} \frac{1}{p} + \sum_{p \leq x} \sum_{n=2}^\infty \frac{1}{np^n} \qquad (18)$$
$$= \log\log x + B_1 + R(x) + \sum_{p \leq x} \sum_{n=2}^\infty \frac{1}{np^n},$$

where $B_1 > 0$ is constant and $R(x)$ is the error term in Theorem 2. Substituting the formula $B_1 = \gamma - \sum_{p \geq 2} \sum_{n=2}^\infty \frac{1}{np^n}$, see [HW, p. 351], into the previous expression and simplifying it return

$$\sum_{p \leq x} \log(1 - 1/p)^{-1} = \log\log x + \gamma + R(x) - \sum_{p > x} \sum_{n=2}^\infty \frac{1}{np^n}. \qquad (19)$$

Since the tail of the power series has the order of magnitude $\sum_{p > x} \sum_{n \geq 2} (np^n)^{-1} = O(1/x)$, the total error term in (19) has the form

$$R(x) - \sum_{p > x} \sum_{n=2}^\infty \frac{1}{np^n} = \begin{cases} O(e^{-.2098(\log x)^{3/5}(\log\log x)^{-1/5}}) & \text{if } x \notin E, \\ \Omega_\pm(x^{-1/2} \log\log\log x / \log x) & \text{if } x \in E, \end{cases} \qquad (20)$$

as $x \to \infty$. Now replacing this into the penultimate expression yields

$$\sum_{p \leq x} \log(1 - 1/p)^{-1} = \begin{cases} \log\log x + \gamma + O(e^{-.2098(\log x)^{3/5}(\log\log x)^{-1/5}}) & \text{if } x \notin E, \\ \log\log x + \gamma + \Omega_\pm(x^{-1/2} \log\log\log x / \log x) & \text{if } x \in E, \end{cases} \qquad (21)$$

as $x \to \infty$. Reversing the logarithm confirms the claim (i). The verifications of (ii), (iii), and (iv) are similar mutatis mutandis. ∎

The established estimates of these products are of the forms $\prod_{p \leq x}(1 \pm 1/p)^{\pm 1} = (c_{\pm} \log x)^{\pm 1} + O(1/\log^2 x)$, and conditional on the Riemann hypothesis, the products are $\prod_{p \leq x}(1 \pm 1/p)^{\pm 1} = (c_{\pm} \log x)^{\pm 1} + O(x^{-1/2} \log x)$, where the constant $c_{\pm} = e^{\gamma}$ or $6e^{\gamma}\pi^{-2}$ see [RS].

The limits of the products $P(x) = (\log x)^{-1} \prod_{p \leq x}(1 - 1/p)^{-1}$ and $Q(x) = (\log x)^{-1} \prod_{p \leq x}(1 + 1/p)^{-1}$ are easy to derive. Basically, the first limit gives the rate of growth of the partial product $\zeta_x(s) = \prod_{p \leq x}(1 - 1/p^s)^{-1}$ of the zeta function $\zeta(s) = \sum_{n \geq 1} n^{-s} = \prod_{p \geq 2}(1 - 1/p^s)^{-1}$ evaluated at $s = 1$ as $x \to \infty$.

***Theorem 7.*** (Martens) The following asymptotic formulas hold:

(i) $\lim_{n \to \infty}(\log n)^{-1} \prod_{p \leq p_n}(1 - 1/p)^{-1} = e^{\gamma}$,

(ii) $\lim_{n \to \infty}(\log n)^{-1} \prod_{p \leq p_n}(1 + 1/p)^{-1} = 6e^{\gamma}/\pi^2$.

The proofs of these results appear in [HW]. The more general products for primes in an arithmetic progressions $qn + a$ with $\gcd(a, q) = 1$, $n \in \mathbb{N}$, have the forms $\prod_{p \leq x, p \equiv a \bmod q}(1 \pm 1/p)^{\pm 1} = (c_{\pm} \log x)^{\pm 1/\varphi(q)} + O_{a,q}(1/\log^2 x)$. These are considered in [LC] and earlier works.

This Section concludes with an application of the product in Theorem 6. This application is concerned with the dynamic of the product $\prod_{p \leq x}(1 - 1/p)^{-1}$ for $x$ close to $e^{\gamma} \log N_k$. This idea might be relevant in estimating the size of the exceptional subset of integers $E \cap \mathbb{N}$.

Let $N_k = 2^{v_1} \cdot 3^{v_2} \cdots p_k^{v_k} \in \mathbb{N}$ be a large primorial integer such that $p_k \leq \beta \log N_k$, where $\beta = \beta(N_k)$ is a parameter and consider the product

$$\prod_{p \leq \beta \log N_k}(1 - 1/p)^{-1}. \tag{22}$$

The oscillation of this product for $x = \beta \log N_k$ near the real number $e^{\gamma} \log N_k$ appears to be a direct function of the parameter $\beta = \beta(N_k) > 0$. Ultra smooth integers with $p_k \leq \beta \log N_k < \log N_k / \log \log N_k$ or similar, seem to have negative oscillations, but ultra smooth integers with $p_k \leq \beta \log N_k < \log N_k \log \log N_k$ or similar, seem to have positive oscillations. The special case $p_k \leq \beta \log N_k \leq c_0 \log N_k$, where $c_0 \approx 1$ is a small constant, known as *colossally abundant* integers $N_k = 2^{v_1} \cdot 3^{v_2} \cdots p_k^{v_k}$ with $v_1 \geq v_2 \geq \cdots \geq v_k \geq 1$, have been investigated extensively, confer the literature for more details.

***Corollary 8.*** Let $N_k = 2^{v_1} \cdot 3^{v_2} \cdots p_k^{v_k} \in \mathbb{N}$ be a large primorial integer such that $p_k \leq c_0 \log N_k$, with $c_0 > 0$ constant and $v_i \geq 1$ for all $1 \leq i \leq k$. Then $N_k / \varphi(N_k) > e^{\gamma} \log \log N_k$.

***Proof***: By hypothesis and Theorem 6, the normalized totient function $N_k / \varphi(N_k) = \prod_{p \mid N_k}(1 - 1/p)^{-1}$ satisfies the formula

$$\prod_{p \leq c_0 \log N_k}(1-1/p)^{-1} = \begin{cases} e^\gamma \log\log c_0 N_k + O(e^{-.2098(\log\log c_0 N_k)^{3/5}(\log\log\log c_0 N_k)^{-1/5}}) & \text{if } x \notin E, \\ e^\gamma \log\log c_0 N_k + \Omega_\pm((\log c_0 N_k)^{-1/2} \log\log\log\log c_0 N_k / \log\log c_0 N_k) & \text{if } x \in E, \end{cases} \quad (23)$$

as $N_k \to \infty$, notice that the effect of the constant $c_0 > 0$ is asymptotically insignificant. Now it is known that colossally abundant integers $N_k$ maximize the divisor function

$$\begin{aligned}\frac{\sigma(N_k)}{N_k} &= \frac{N_k}{\varphi(N_k)} \prod_{p^\alpha \| N_k}(1-1/p^{\alpha+1}) \\ &= \prod_{p \leq c_0 \log N_k}(1-1/p)^{-1} \prod_{p^\alpha \| N_k}(1-1/p^{\alpha+1}).\end{aligned} \quad (24)$$

Ergo, it readily follows that the extreme values

$$\begin{aligned}N_k/\varphi(N_k) &= e^\gamma \log\log N_k + \Omega_\pm\left((\log N_k)^{-1/2} \log\log\log\log N_k / \log\log N_k\right) \\ &\geq e^\gamma \log\log N_k + c_1(\log N_k)^{-1/2} \log\log\log\log N_k / \log\log N_k \\ &> e^\gamma \log\log N_k\end{aligned} \quad (25)$$

hold for infinitely many integers $N_k$, with $c_1 > 0$ constant. ∎

For squarefree integers $N_k$, the sum of divisor function in (24) reduces to

$$\begin{aligned}\frac{\sigma(N_k)}{N_k} &= \frac{N_k}{\varphi(N_k)} \prod_{p \| N_k}(1-1/p^2) \\ &> \frac{6}{\pi^2} \prod_{p \leq c_0 \log N_k}(1-1/p)^{-1} \\ &> \frac{6e^\gamma}{\pi^2} \log\log N_k.\end{aligned} \quad (26)$$

These results show that the extreme values of the product $\prod_{p \leq x}(1-1/p)^{-1}$ occur at the colossally abundant integers. $N_k = 2^{v_1} \cdot 3^{v_2} \cdots p_k^{v_k}$ with $v_1 \geq v_2 \geq \cdots \geq v_k \geq 1$, and possibly other ultra smooth integers, see [NS] for related works. So it is possible that the exceptional subset of integers $E \cap \mathbb{N} = \{$ colossally abundant integers $\}$?

**4. The Proof of Theorem 1**
The proof of Theorem 1 is obtained from Theorem 6. To show this, observe that the oscillatory component of the product is

$$\prod_{p \leq x}(1-1/p)^{-1} = e^\gamma \log x + \Omega_\pm(x^{-1/2} \log\log\log x / \log x) \quad (27)$$

for $x \in E$, see Definition 5. Therefore,

$$x^{1/2} \log x \left( \prod_{p \leq x} (1 - 1/p)^{-1} - e^{\gamma} \log x \right) = \Omega_{\pm}(\log \log \log x). \tag{28}$$

Clearly, this implies that the product $\prod_{p \leq x}(1 - 1/p)^{-1}$ oscillates symmetrically about $e^{\gamma} \log x$ as $x \to \infty$.

## 5 More Finite Sums Over The Prime Numbers

The same technique employed in Section 3 will be used to derive estimates of some common finite sums over the prime numbers. The earliest finite sums over the prime numbers of interest are the theta function $\vartheta(x) = \sum_{p \leq x} \log p$ and the psi function $\psi(x) = \sum_{n \leq x} \Lambda(n)$ for $x \geq 1$, where the number theoretical function is defined by

$$\Lambda(n) = \begin{cases} \log p & \text{if } n = p^{\nu}, \\ 0 & \text{if } n \neq p^{\nu}, \end{cases} \tag{29}$$

for prime powers $n = p^{\nu}$, $\nu \geq 1$. These two functions are linked via the Mobius inverse identities $\psi(x) = \sum_{n \geq 1} \vartheta(x^{1/n})$ and $\vartheta(x) = \sum_{n \geq 1} \mu(n) \psi(x^{1/n})$. The Mobius pair leads to the asymptotic formula $\psi(x) = \vartheta(x) + O(x^{1/2})$.

**Theorem 9.** Let $x > 1$ be a large number. Then

$$\vartheta(x) = \begin{cases} x + O(xe^{-.2098(\log x)^{3/5}(\log \log x)^{-1/5}}) & \text{if } x \notin E, \\ x + \Omega_{\pm}(x^{1/2} \log \log \log x) & \text{if } x \in E. \end{cases} \tag{30}$$

**Proof:** Use the definition of the theta function $\vartheta(x) = \sum_{p \leq x} \log p$, and integration by parts of the integral representation

$$\sum_{p \leq x} \log p = \sum_{n \leq x} (\pi(n) - \pi(n-1)) \log n = \int_2^x (\log t) d\pi(t), \tag{31}$$

to prove the claim. ∎

A discussion of the earliest works on theta and psi functions appears in [IV, p. 347], and [SP, p. 451].

**Theorem 10.** Let $x > 1$ be a large number, and let $\alpha \in \mathbb{R}$ be a real number. Then the twisted theta function $\vartheta(\alpha, x) = \sum_{p \leq x} e^{i\alpha p} \log p$ satisfies the following:

$$\vartheta(\alpha, x) = \begin{cases} e^{i\alpha x} x + O(xe^{-.2098(\log x)^{3/5}(\log \log x)^{-1/5}}) & \text{if } x \notin E, \\ e^{i\alpha x} x + \Omega_{\pm}(x^{1/2} \log \log \log x) & \text{if } x \in E. \end{cases} \tag{32}$$

**Proof:** Utilize integration by parts of the integral representation

$$\sum_{p \leq x} e^{i\alpha p} \log p = \sum_{n \leq x} (\pi(n) - \pi(n-1))e^{i\alpha n} \log n = \int_2^x (\log t)e^{i\alpha t} d\pi(t), \tag{33}$$

to prove the claim. ∎

This technique easily shows that the finite exponential sum

$$\sum_{p \leq x} e^{i\alpha p} = \begin{cases} e^{i\alpha x} \pi(x) + O(xe^{-.2098(\log x)^{3/5}(\log\log x)^{-1/5}}) & \text{if } x \notin E, \\ e^{i\alpha x} \pi(x) + \Omega_{\pm}(x^{1/2} \log\log\log x / \log x) & \text{if } x \in E, \end{cases} \tag{34}$$

attains large positive and negative values infinitely often, and provides explicit asymptotic evaluations. For example, at $\alpha = \pi/2$, and $x = 2n + 1$, $n \geq 0$, the finite sum takes the shape

$$\sum_{p \leq x} e^{i\alpha p} = \begin{cases} \pm i\pi(x) + O(xe^{-.2098(\log x)^{3/5}(\log\log x)^{-1/5}}) & \text{if } x \notin E, \\ \pm i\pi(x) + \Omega_{\pm}(x^{1/2} \log\log\log x / \log x) & \text{if } x \in E. \end{cases} \tag{35}$$

**Theorem 11.** Let $x > 1$ be a large number. Then

$$\sum_{p \leq x} \frac{\log p}{p} = \begin{cases} \log x - c_0 + O(e^{-.2098(\log x)^{3/5}(\log\log x)^{-1/5}}) & \text{if } x \notin E, \\ \log x - c_0 + \Omega_{\pm}(x^{-1/2} \log\log\log x) & \text{if } x \in E, \end{cases} \tag{36}$$

where $c_0 > 0$ is a constant.

*Proof*: Use the previous Theorem to compute the integral representation

$$\sum_{p \leq x} \frac{\log p}{p} = \sum_{n \leq x} \frac{(\vartheta(n) - \vartheta(n-1))}{n} = \int_2^x \frac{d\vartheta(t)}{t}. \tag{37}$$

This completes the verification. ∎

This result is usually stated as $\sum_{p \leq x} \log p / p = \log x + O(1)$ without any information on the error term $O(1)$, see the literature. The same estimate (35) also holds for the closely related finite sum $\sum_{n \leq x} \Lambda(n)/n$, where $\Lambda(n) = \log p$ if $n = p^v \geq 2$ is a prime power, else $\Lambda(n) = 0$. This finite sum coincides with the partial sum of the logarithm of the zeta function evaluated at $s = 1$, viz, $-\zeta'(s)/\zeta(s) = \sum_{n \geq 1} \Lambda(n)/n^{-s}$.

**Theorem 12.** Let $x > 1$ be a large number, and let $((x))$ be the fractional part of $x$. Then

$$\sum_{p \leq x} ((x/p)) \log p = \begin{cases} (1 - c_0)x + O(xe^{-.2098(\log x)^{3/5}(\log\log x)^{-1/5}}) & \text{if } x \notin E, \\ (1 - c_0)x + \Omega_{\pm}(x^{1/2} \log\log\log x) & \text{if } x \in E. \end{cases} \tag{38}$$

***Proof***: The idea is to compute the logarithm of the factorial function $x!$ in two ways to derive the result. Consider the prime decomposition $x! = \prod_{p \leq x} p^{v_p}$, where $v_p = \sum_{n \geq 1} [x/p^n]$, and $[\,x\,]$ is the largest integer function. The first computation is

$$\log \prod_{p \leq x} p^{v_p} = \sum_{n \geq 1} \sum_{p \leq x} [x/p^n] \log p$$
$$= \sum_{p \leq x} [x/p] \log p + \sum_{n \geq 2} \sum_{p \leq x} [x/p^n] \log p \qquad (39)$$
$$= x \sum_{p \leq x} \frac{\log p}{p} - \sum_{p \leq x} ((x/p)) \log p + O(\log^2 x).$$

The last line uses $[\,x\,] = x - ((\,x\,))$ and the estimate $\sum_{n \geq 2} \sum_{p \leq x} [x/p^n] \log p = O(\log^2 x)$. Now applying Theorem 10 yields

$$\log \prod_{p \leq x} p^{v_p} = x(\log x - c_0) + E(x) - \sum_{p \leq x} ((x/p)) \log p + O(\log^2 x), \qquad (40)$$

where

$$E(x) = \begin{cases} O(xe^{-.2098(\log x)^{3/5}(\log\log x)^{-1/5}}) & \text{if } x \notin E, \\ \Omega_{\pm}(x^{1/2} \log\log\log x) & \text{if } x \in E. \end{cases} \qquad (41)$$

.
The second computation is

$$\log x! = (x + 1/2) \log x - x + c_1 + O(1/x), \qquad (42)$$

where $c_0, c_1 > 0$ are constants. Matching (39) and (41) and regrouping the terms return the claim. ∎

The constant $c_0$ should be in the range $0 \leq c_0 \leq \log 2$. Surprisingly, the finite sums $\sum_{p \leq x} \log p$ and $\sum_{p \leq x} ((x/p)) \log p$ have almost the same asymptotic behavior. Moreover, note that for an integer $x \in \mathbb{N}$, the equation $((x/p)) = 0$ has exactly $\omega(x) \geq 1$ zeros, otherwise $1/p \leq ((x/p)) \leq (p-1)/p$, so the finite sum has a nonzero minimal $\sum_{p \leq x}((x/p)) \log p \geq \sum_{p \leq x} (\log p)/p - \omega(x) \geq \log x - O(\log x / \log\log x)$. Here the omega symbol is defined by $\omega(n) = \#\{\, p \mid n : p \text{ prime }\,\}$, which is the number of prime divisors function. Moreover, assuming the Riemann Hypothesis, this estimate can be sharpened to

$$\sum_{p \leq x}((x/p)) \log p = \begin{cases} (1 - c_0)x + O(x^{1/2} \log x) & \text{if } x \notin E, \\ (1 - c_0)x + \Omega_{\pm}(x^{1/2} \log\log\log x) & \text{if } x \in E. \end{cases} \qquad (43)$$

Let $s \in \mathbb{C}$ be a complex number. The prime harmonic sum $\sum_{p \leq x} p^{-s}$ generalizes the special case $\sum_{p \leq x} 1/p$. The limit $\zeta^*(s) = \sum_{p \geq 2} p^{-s}$ is called the prime zeta function. For computational works, the faster convergent formula $\zeta^*(s) = \sum_{n \geq 1} \mu(n) n^{-s} \log \zeta(ns)$ can be used.

***Theorem* 13.**  Let $x \geq x_0$, and let $s \in \mathbb{C}$ be a fixed complex number, $\mathfrak{Re}(s) \neq 1$. Then

$$\sum_{p \leq x} \frac{1}{p^s} = \begin{cases} x^{-s-1} \log\log x + c_1 + O(x^{1-s} e^{-.2098(\log x)^{3/5}(\log\log x)^{-1/5}}) & \text{if } x \notin E, \\ x^{-s-1} \log\log x + c_0 + \Omega_{\pm}(x^{-1/2-s} \log\log\log x / \log x) & \text{if } x \in E, \end{cases} \quad (44)$$

where $c_1 > 0$ is a constant which depends on $s$, and $E \subset \mathbb{R}$ is the exceptional subset of real numbers.

***Proof*:** The integral representation of the prime harmonic sum is

$$\sum_{p \leq x} \frac{1}{p^s} = \sum_{n \leq x} \frac{\pi(n) - \pi(n-1)}{n^s} = \int_2^x \frac{d\pi(t)}{t^s}. \quad (45)$$

Using the exceptional part of the prime counting measure (8) yields

$$\begin{aligned}
\sum_{p \leq x} \frac{1}{p^s} &= \int_2^x \frac{d(li(t) + \Omega_{\pm}(t^{1/2} \log\log\log t / \log t))}{t^s} \\
&= \int_2^x \frac{dt}{t^s \log t} + \int_2^x \frac{d\Omega_{\pm}(t^{1/2} \log\log\log t / \log t)}{t^s} \\
&= x^{-s+1} \log\log x + c_3 + \frac{\Omega_{\pm}(x^{1/2} \log\log\log x / \log x)}{x^s} - \int_2^x \frac{\Omega_{\pm}(t^{1/2} \log\log\log t / \log t) dt}{t^{s+2}},
\end{aligned} \quad (46)$$

where $c_2 = -2^{-s+1}\log\log 2$, $c_3 > 0$ are constants. The last integral is written in the form

$$\int_2^x \frac{\Omega_{\pm}(t^{1/2} \log\log\log t / \log t) dt}{t^{2+s}} = \int_2^\infty \frac{\Omega_{\pm}(t^{1/2} \log\log\log t / \log t) dt}{t^{2+s}} - \int_x^\infty \frac{\Omega_{\pm}(t^{1/2} \log\log\log t / \log t) dt}{t^{2+s}} \\ = c_4 + \Omega_{\pm}(x^{-1/2-s} \log\log\log x / \log x). \quad (47)$$

The claim follows from these observations. ∎

Given a fixed complex number $1 \neq s \in \mathbb{C}$, the prime harmonic series either converges to a constant $c_2 = c_2(s)$ or it diverges depending on the real part $\mathfrak{Re}(s)$.